\documentclass{ACmod}

\DeclareFontFamily{U}{rsf}{}
\DeclareFontShape{U}{rsf}{m}{n}{
  <5> <6> rsfs5 <7> <8> <9> rsfs7 <10-> rsfs10}{}
\DeclareMathAlphabet{\mathscr}{U}{rsf}{m}{n}

\DeclareMathAlphabet{\mathgth}{U}{euf}{m}{n}

\DeclareFontFamily{U}{cyr}{}
\DeclareFontShape{U}{cyr}{m}{n}{
  <5> wncyr5 <6> wncyr6 <7> wncyr7 <8> wncyr8 <9> wncyr9 <10-> wncyr10}{}
\DeclareMathAlphabet{\mathcyr}{U}{cyr}{m}{n}

\input cyracc.def

\DeclareSymbolFont{bbold}{U}{bbold}{m}{n}
\DeclareSymbolFontAlphabet{\mathbbold}{bbold}

\makeatletter
\def\operator@font{\sf}
\makeatother

\renewcommand{\S}{{\mathsf{Sym}}}

\renewcommand{\phi}{\varphi}

\usepackage{fancyhdr}
\usepackage{amssymb}
\usepackage{xypic}

\usepackage{tikz-cd}
\usepackage{mathtools}
\usepackage{amsfonts}

\usepackage{tikz}
\usepackage[tocflat]{tocstyle}
\usetikzlibrary{positioning}
\usepackage{amscd}

\usepackage{amsthm}
\usepackage{amsmath}

\everymath{\displaystyle}

\begin{document}

\title{Functoriality of HKR isomorphisms}

\author[Shengyuan Huang]{%
Shengyuan Huang}

\address{Shengyuan Huang, Mathematics Department,
University of Wisconsin--Madison, 480 Lincoln Drive, Madison, WI
53706--1388, USA.}

\begin{abstract}
  {\sc Abstract:} For a closed embedding of smooth schemes $X\hookrightarrow S$ with a fixed first order splitting, one can construct HKR isomorphisms between the derived scheme $X\times^R_S X$ and the total space of the shifted normal bundle $\mathbb{N}_{X/S}[-1]$, due to Arinkin-C\u{a}ld\u{a}raru, Arinkin-C\u{a}ld\u{a}raru-Hablicsek, and Grivaux. In this paper, we study functoriality property of the HKR isomorphisms for a sequence of closed embeddings $X\hookrightarrow Y\hookrightarrow S$.  The HKR isomorphism is functorial when a certain cohomology class, which we call the Bass-Quillen class, vanishes.  We obtain Lie theoretic interpretations for the HKR isomorphisms and for the Bass-Quillen class as well.
\end{abstract}

\maketitle

\setcounter{tocdepth}{1}
\tableofcontents

\section{Introduction}

\paragraph
Let $X$ be a smooth algebraic variety over a field of characteristic zero.  There is a Hochschild-Kostant-Rosenberg isomorphism \cite{S} in the derived category of $X$
$$\Delta^*\Delta_*\mathcal{O}_X\cong\S_{\mathcal{O}_X}(\Omega_X[1]),$$
where $\Delta$ is the diagonal embedding of $X$ in $X\times X$.

It was observed by Kapranov and Kontsevich that there is a Lie theoretic interpretation of the HKR isomorphism.  The derived loop space $LX=X\times^R_{X\times X}X$ has the structure of a derived group scheme over $X$ and the shifted tangent bundle $T_X[-1]$ is its Lie algebra \cite{K}.  The HKR isomorphism can be thought of as a version of the exponential map $\mathbb{T}_X[-1]\rightarrow LX$  \cite{CR}, where $\mathbb{T}_X[-1]$ is the total space of the shifted tangent bundle.
\paragraph In the case of classical Lie groups we have a commutative diagram
$$
\xymatrix{
G\ar[r]^{f} & H\\
\mathfrak{g}\ar[u]^{\mathrm{exp}}\ar[r]^{df} & \mathfrak{h}\ar[u]_{\mathrm{exp}}. }
$$
A map of schemes $X\rightarrow Y$ induces a map of derived group schemes $LX\rightarrow LY|_X$ over $X$.  The analogous statement for derived schemes to the above Lie theoretic statement is that the diagram
$$
\xymatrix{
LX\ar[r] & LY|_X\\
\mathbb{T}_X[-1]\ar[r]\ar[u]^{\mbox{HKR}} & \mathbb{T}_Y|_X[-1]\ar[u]_{\mbox{HKR}}
}
$$
commutes.  One can prove the commutativity of the diagram above easily using the methods in this paper.
\paragraph We would like to consider the commutativity of this diagram in a more general setting.  One can still get HKR isomorphisms if one replaces the diagonal embedding by an arbitrary closed embedding of schemes $i:X\hookrightarrow S$ with a fixed first order splitting. Arinkin, C\u{a}ld\u{a}raru, Hablicsek \cite{AC, ACH}, and Grivaux \cite{Gri} constructed HKR isomorphisms between the total space of the shifted normal bundle $\mathbb{N}_{X/S}[-1]$ and the derived self-intersection $X\times^R_SX$ from a fixed first order splitting of $i$.  In this paper we study the functoriality of the HKR isomorphism constructed in \cite{AC}.

\paragraph We will be in the following setting from now on.  Let $X\hookrightarrow Y\hookrightarrow S$ be a sequence of closed embeddings of smooth schemes.  Assume that $X$ is split to first order in $Y$, and similarly for $X$ in $Y$ and $Y$ in $S$.  We want to understand if the two diagrams
$$
\xymatrix{
X\times^R_YX\ar[r] & X\times^R_SX\ar[r] & Y\times^R_SX=Y\times^R_SY|_X  \\
\mathbb{N}_{X/Y}[-1]\ar[r]\ar[u]^{\cong} & \mathbb{N}_{X/S}[-1]\ar[r]\ar[u]^{\cong} & \mathbb{N}_{Y/S}|_X[-1]\ar[u]^{\cong}
}
$$
are commutative.

\paragraph There are two ways to construct isomorphisms between $\mathbb{N}_{X/S}[-1]$ and $X\times^R_SX$ from a fixed first order splitting.  The construction of HKR isomorphisms in \cite{AC} involves a differential graded $\mathcal{O}_{X^{(1)}_S}$-algebra resolution of $\mathcal{O}_X$.  The construction of $\mbox{HKR}$ isomorphisms in \cite{ACH} and \cite{Gri} reduces all HKR isomorphisms to $\mbox{HKR}$ isomorphisms of diagonal embeddings constructed in \cite{AC}.  For more details see section 2.  It is unknown whether the first and the second constructions agree or not.  In this paper we only deal with the HKR isomorphisms constructed in \cite{AC}.

\paragraph The following question is important in this paper.  The restriction of the conormal bundle $N^\vee_{Y/S}$ to the first order neighborhood $X^{(1)}_Y$ is a vector bundle on $X^{(1)}_Y$.  One can ask whether this bundle is isomorphic to $s^*(N^\vee_{Y/S}|_X)$, where $s$ is the chosen first order splitting of $X$ in $Y$.  If the answer to this question is positive, i.e., there is an isomorphism of vector bundles on $X^{(1)}_Y$
      $$N^\vee_{Y/S}|_{X^{(1)}_Y}\cong s^*(N^\vee_{Y/S}|_X),\mbox{ }\mbox{ }(*)$$
we will say condition $(*)$ is satisfied.
\paragraph
\label{subsec:BQ} The condition $(*)$ is equivalent to the vanishing of a certain cohomology class in $\mathrm{Ext}^1(N_{X/Y}\otimes N_{Y/S}|_X, N_{Y/S}|_X)$ associated to $N_{Y/S}^\vee|_{X^{(1)}_Y}$.  One can construct a similar cohomology class $\alpha_{s,\mathcal{G}}\in\mathrm{Ext}^1(N_{X/Y}\otimes\mathcal{G}^\vee|_{X}, \mathcal{G}^\vee|_{X})$ for any vector bundle $\mathcal{G}$ on $X^{(1)}_Y$ and the fixed splitting $s$.  This class vanishes if and only if $\mathcal{G}\cong s^*(\mathcal{G}|_X)$.  We will call this cohomology class the Bass-Quillen class since it is related to the Bass-Quillen conjecture as we explain below.

Suppose $X$ is a smooth algebraic variety, and $Y$ is the total space of a vector bundle $\mathcal{E}$ on $X$.  Let $\mathcal{F}$ be a vector bundle on $Y$.  Then one can ask if $\mathcal{F}$ is isomorphic to the pull back of some vector bundle on $X$.  When $X$ is affine this is known as the Bass-Quillen problem and was answered affirmatively in \cite{L}.  However, in the global case the answer to this question is negative.  In particular, there can be no vector bundle on $X$ whose pull back is $\mathcal{F}$ if the Bass-Quillen class in $\mathrm{Ext}^1(N_{X/Y}\otimes\mathcal{F}^\vee|_{X},\mathcal{F}^\vee|_{X})$ associated to $\mathcal{F}|_{X^{(1)}_Y}$ is not zero.

Here are the main results in this paper.

\paragraph{\bf Theorem A.}
{\em
Let $X\hookrightarrow Y\hookrightarrow S$ be a sequence of closed embeddings of smooth schemes.  Further assume that there are compatible splittings on the tangent bundles
$$
\xymatrix{
T_X \ar[r] & T_Y|_X \ar[r]\ar@{-->}@/_1.0pc/[l]_{p} & T_S|_X. \ar@{-->}@/_1.0pc/[l]_{q|_X}\ar@{-->}@/^1.0pc/[ll]^{\rho}
}
$$
Compatibility means that $p\circ q|_X=\rho$.

 The left square below is commutative.  If condition $(*)$ is satisfied, i.e., the Bass-Quillen class of $N^\vee_{Y/S}|_{X^{(1)}_Y}$ is zero, then the right square below is also commutative.
$$
\xymatrix{
X\times^R_YX\ar[r] & X\times^R_SX\ar[r] & Y\times^R_SX=Y\times^R_SY|_X  \\
\mathbb{N}_{X/Y}[-1]\ar[r]\ar[u]^{\cong} & \mathbb{N}_{X/S}[-1]\ar[r]\ar[u]^{\cong} & \mathbb{N}_{Y/S}|_X[-1]\ar[u]^{\cong}.
}
$$
The vertical maps are the HKR isomorphisms defined in \cite{AC}.  The horizontal maps between normal bundles are the linear ones, i.e., they are vector bundle maps.
   }

\paragraph{\bf Application to orbifolds.} This paper is motivated by the study of Hochschild cohomology of an orbifold. The main example we consider is the setting where $S$ admits an action of a finite group $G$ and $X$ is the fixed locus of $G$ and $Y$ is the fixed locus of a subgroup $H\leq G$.  It is easy to see that the splittings obtained from the averaging maps are compatible in this case.

\paragraph There is a Lie theoretic interpretation for Theorem A.  Under certain extra assumptions, all derived self-intersections in Theorem A are groups in the derived category of dg schemes.  The shifted normal bundles are their Lie algebras.  One can check that the natural maps
$$X\times^R_YX\longrightarrow X\times^R_SX\longrightarrow X\times^R_SY=Y\times^R_SY|_X$$
are maps of groups.  The map $N_{X/Y}[-1]\rightarrow N_{Y/S}|_X[-1]$ respects the Lie structures in general.  However, $N_{X/S}[-1]\rightarrow N_{Y/S}|_X[-1]$ may not preserve the Lie brackets.  The Bass-Quillen class $\alpha_{s,N^\vee_{Y/S}|_{X^{(1)}_Y}}$ is precisely the obstruction to this map preserving the Lie brackets.

With no assumptions other than the choice of a first order splitting, the bundle $N_{X/S}[-1]$ carries an anti-symmetric bracket which may not satisfy the Jacobi identity.  We call this structure a pre-Lie bracket.  We will provide more details in section 2.

{\bf Theorem B.} In the same setting as Theorem A, the vector bundle map $N_{X/S}[-1]\rightarrow N_{Y/S}|_X[-1]$ preserves the pre-Lie brackets if and only if the Bass-Quillen class $\alpha_{s,N^\vee_{Y/S}|_{X^{(1)}_Y}}:N_{X/Y}\otimes N_{Y/S}|_X\rightarrow N_{Y/S}|_X[1]$ is zero.

Thus Theorem A and Theorem B provide a generalization of the original result for Lie groups to the setting of groups obtained as self-intersections.

\paragraph {\bf Plan of the paper.} Section 2 is reviewing background material.  We recall what is known about the structure of arbitrary closed embeddings with further assumptions on (possibly higher order) splittings, and how the $\mbox{HKR}$ isomorphisms are constructed for a closed embedding with a fixed first order splitting.  Then we define the Bass-Quillen class.  We also provide details about the Lie theoretic interpretation of the Bass-Quillen class and our Theorem B.

Section 3 is devoted to the proof of theorem A.  To prove commutativity of the diagram in Theorem A, it suffices to check whether the explicit resolutions of $\mathcal{O}_X$ and $\mathcal{O}_Y$ in \cite{AC} are compatible.

Section 4 is about Theorem B for Lie algebras in the category of vector spaces.  We consider an inclusion of Lie algebras $\mathfrak{h}\hookrightarrow\mathfrak{g}$ and the associated short exact sequence of $\mathfrak{h}$-modules
$$
\xymatrix{
 0\ar[r] & \mathfrak{h} \ar[r]^{\alpha} & \mathfrak{g} \ar[r]^{\beta} & \mathfrak{n}=\mathfrak{g}/\mathfrak{h}\ar[r]&0.
}
$$
There is a way to construct a pre-Lie bracket on $\mathfrak{n}$ which may not be compatible with the Lie bracket on $\mathfrak{g}$.  Then we prove a result similar to Theorem B using methods that make sense in the derived category. This approach will be generalized to derived category of schemes in section 6.

Section 5 is the proof of our Lie theoretic interpretation of the Bass-Quillen class $\alpha_{s,N^\vee_{Y/S}|_{X^{(1)}_Y}} $.  This class can be viewed as a Lie module structure map, where $N_{X/Y}[-1]$ is the Lie algebra and $N_{Y/S}|_X[-1]$ is the module.

Section 6 generalizes the proofs in section 4 to the derived setting.  We end this paper with an example where the Bass-Quillen class is not zero.

\paragraph {\bf Conventions.} All the schemes we considered in this paper are smooth over a field of characteristic zero.

\paragraph {\bf Acknowledgments.} I would like to thank and Andrei C\u{a}ld\u{a}raru for discussing details with me during our weekly meeting.  I am also grateful to Dima Arinkin for patiently listening to the various problems we ran into at different stages of the project, and for providing insight.

\section{Background and Lie theoretic interpretations}

We first discuss what is known about the HKR isomorphism for the diagonal embedding.  Then we recall the definition of the HKR isomorphisms for a general closed embedding $X\hookrightarrow S$ with a fixed first order splitting.  We briefly recall the Lie theoretic interpretations of general embeddings with possibly higher order splittings.  We define the Bass-Quillen class and explain the Lie theoretic interpretations for the Bass-Quillen class $\alpha_{s,N^\vee_{Y/S}|_{X^{(1)}_Y}}$ and Theorem B at last.

\paragraph{\bf The diagonal embedding.} Let $X$ be a smooth algebraic variety.  There is an HKR isomorphism
$$\mathrm{HH}^*(X)\cong\bigoplus_{p+q=*}H^{p}(X,\wedge^q T_X)$$
that identifies Hochschild cohomology of $X$ with polyvector fields as vector spaces.  More precisely, we have an HKR isomorphism at the level of sheaves in the derived category of $X$
$$\Delta^*\Delta_*\mathcal{O}_X\cong\S_{\mathcal{O}_X}(\Omega_X[1]),$$
where $\Delta: X\hookrightarrow X\times X$ is the diagonal embedding.  We get the desired isomorphism on cohomology by applying $\mathrm{RHom}(-,\mathcal{O}_X)$ to the isomorphism of sheaves.
\paragraph In the world of derived schemes we consider the free loop space $LX$ of $X$, defined as the derived self-intersection $X\times^R_{X\times X} X$.  Its structure complex is $\Delta^*\Delta_*\mathcal{O}_X$, and the structure complex of the total space of the shifted tangent bundle $\mathbb{T}_{X}[-1]=\mathrm{Spec}_{\mathcal{O}_X}(\S(\Omega_X[1]))$ is $\S(\Omega_X[1])$.  We can restate the HKR isomorphism as an isomorphism of derived schemes over $X$
$$\xymatrix{
\mathbb{T}_X[-1]\ar[r]^{\cong~~~~~~} & LX=X\times^R_{X\times X} X.
}
$$
It can be viewed as the exponential map from the Lie algebra $\mathbb{T}_{X}[-1]$ to the group $LX$ as explained in the introduction.
\paragraph{\bf General embeddings.} There exist generalized HKR isomorphisms if we replace the diagonal embedding by an arbitrary closed embedding $i: X\hookrightarrow S$ of smooth schemes.

The embedding $i$ factors as
$$
\xymatrix{
X\ar@{^{(}->}[r]^{\mu} & X^{(1)}_S\ar@{^{(}->}[r]^{\nu} & S,
}
$$
where $X^{(1)}_S$ is the first order neighborhood of $X$ in $S$.  We say $i$ splits to first order if and only if the map $\mu$ is split, i.e., there exists a map of schemes $\varphi: X^{(1)}_S\rightarrow X$ such that $\varphi\circ\mu=id$.  There is a bijection between first order splittings of $i$ and splittings of the short exact sequence below \cite[20.5.12 (iv)]{G}
$$
\xymatrix{
 0\ar[r] & T_X \ar[r] & T_S|_X \ar[r]\ar@{-->}@/_1.0pc/[l] & N_{X/S}\ar[r] & 0.
 }
$$

\paragraph Arinkin and C\u{a}ld\u{a}raru \cite{AC} provided a necessary and sufficient condition for $i^*i_*\mathcal{O}_X$ to be isomorphic to $\S(N_{X/S}^\vee[1])$.  In \cite{ACH} Arinkin, C\u{a}ld\u{a}raru, and Hablicsek proved that the derived intersection $X\times^R_SX$ is isomorphic to $\mathbb{N}_{X/S}[-1]$ over $X\times X$ if and only if the embedding $i$ splits to first order.  Grivaux independently proved a similar result for complex manifolds in \cite {Gri}.

\paragraph\label{subsec: definition of HKR_1}
Let us briefly recall how the $\mbox{HKR}$ isomorphism $i^*i_*\mathcal{O}_X\cong\S(N^\vee_{X/S}[1])$ was constructed in \cite{AC}.  It is defined as the composite map
$$
\xymatrix{
\mu^*\nu^*\nu_*\mu_*\mathcal{O}_X\ar[r] & \mu^*\mu_*\mathcal{O}_X\ar[r]^{\cong} & \mathrm{T}^c(N^\vee_{X/S}[1])\ar[r]^{\exp} & \mathrm{T}(N^\vee_{X/S}[1])\ar[r] &\S(N^\vee_{X/S}[1]).
}
$$
The left most map is given by the counit of the adjunction $\nu^*\dashv\nu_*$.  The map $\exp$ is multiplying by $1/k!$ on the degree $k$ piece, and the last one is the natural projection map.  The $\mathrm{T}^c(N^\vee_{X/S}[1])$ is the free coalgebra on $N^\vee_{X/S}$ with the shuffle product structure, and $\mathrm{T}(N^\vee_{X/S}[1])$ is the tensor algebra on $N^\vee_{X/S}$.  The isomorphism $\mu^*\mu_*\mathcal{O}_X\cong\mathrm{T}^c(N^\vee_{X/S}[1])$ in the middle is non-trivial and needs more explanation.  With the splitting $\varphi$ one can build an explicit resolution of $\mu_*\mathcal{O}_X$ as an $\mathcal{O}_{X^{(1)}_S}$-algebra
$$
 \xymatrix{
 (\mathrm{T}^c(\varphi^*N^\vee_{X/S}[1]),d)\ar[r] & \mu_*\mathcal{O}_{X},
 }
 $$
where $(\mathrm{T}^c(\varphi^*N^\vee_{X/S}[1]),d)$ is the free coalgebra on $\varphi^*N^\vee_{X/S}$ with the shuffle product structure and a differential $d$.  The differential is defined as follows.  There is a short exact sequence on $X^{(1)}_S$
$$0\rightarrow \mu_*N^\vee_{X/S}\rightarrow\mathcal{O}_{X^{(1)}_S}\rightarrow\mu_*\mathcal{O}_X\rightarrow0.$$

Consider the composite map
$$\varphi^*N^\vee_{X/S} \rightarrow \mu_*\mu^*\varphi^*N^\vee_{X/S} =\mu_*N^\vee_{X/S} \rightarrow\mathcal{O}_{X^{(1)}_S},$$
 whose cokernel is $\mu_*\mathcal{O}_X$.  Tensor the morphism above with $(\varphi^*N^\vee_{X/S})^{\otimes (k-1)}$.
 We get the degree $k$-th piece of the differential $d_k: (\varphi^*N^\vee_{X/S})^{\otimes k}\rightarrow(\varphi^*N^\vee_{X/S})^{\otimes (k-1)}$.  The differential vanishes once we pull this resolution back on $X$ via $\mu$, so we get the desired isomorphism.

\paragraph For any vector bundle $\mathcal{E}$ on $X$, we tensor the resolution above by $\varphi^*\mathcal{E}$.  Using the projection formula and $\varphi\circ\mu=id$, one can show that we get a resolution of $\mu_*\mathcal{E}$
$$
(\mathrm{T}^c(\varphi^*N^\vee_{X/S}[1])\otimes\varphi^*\mathcal{E},d)\rightarrow\mu_*\mathcal{E}.
$$
The same argument shows that $i^*i_*(\mathcal{E})\cong\mathcal{E}\otimes\S(N^\vee_{X/S}[1])$, i.e., that $i^*i_*(-)\cong(-)\otimes\S(N^\vee_{X/S}[1])$ as dg functors.  This shows that $X\times^R_S X\cong\mathbb{N}_{X/S}[-1]$ over $X\times X$.

\paragraph \label{subsec: definition of HKR_2}Let us recall how $\mbox{HKR}: \mathbb{N}_{X/S}[-1]\cong X\times^R_SX$ was constructed in \cite{ACH} and \cite{Gri}.  It is defined as the composite map
$$
\xymatrix{\mathbb{N}_{X/S}[-1]\ar@{-->}[r] & \mathbb{T}_S|_X[-1] \ar[r]^{\cong} &  S\times^R_{S\times S}S|_X\ar[dll]^{=} \\ S\times^R_{S\times S}X\ar[r]_{id_S\times\Delta~~~~} & S\times^R_{S\times S}(X\times X)\ar[r]^{\cong} & X\times^R_SX.
}
$$
The dotted arrow is the splitting we fixed.  The isomorphism in the middle $\mathbb{T}_S[-1]\cong S\times^R_{S\times S}S$ is the $\mbox{HKR}$ isomorphism of diagonal embeddings $S\hookrightarrow S\times S$ discussed in~(\ref{subsec: definition of HKR_1}). There are two splittings to define $\mbox{HKR}$ for the diagonal embeddings.  We always choose $p_1$, i.e., the projection onto the left factor
$$\xymatrix{
\Delta_S: S\ar[r] & S\times S\ar@{-->}@/^1.0pc/[l]^{p_2} \ar@{-->}@/_1.0pc/[l]^{p_1}.
}
$$
\paragraph As mentioned in the introduction, we do not know whether the constructions in~(\ref{subsec: definition of HKR_1}) and~(\ref{subsec: definition of HKR_2}) define the same isomorphism or not.  We only consider the one in~(\ref{subsec: definition of HKR_1}) in this paper.

\paragraph{\bf Lie theoretic interpretations for general self-intersections.}
Consider a closed embedding $i:X\hookrightarrow S$ of smooth schemes.  The derived self-intersection $X\times^R_SX$ has an $\infty$-groupoid structure in the $(\infty,1)$-category of dg schemes over $X$.  The associated $L_{\infty}$ algebroid is $N_{X/S}[-1]$.  Passing to the derived category, we get a groupoid in the derived category of dg schemes having $X$ as the space of objects.  The target and source maps are the two projections $\pi_1$ and $\pi_2:X\times^R_SX\rightarrow X$.  See \cite{CCT} for more details.
When $S=X\times X$ and $i$ is the diagonal embedding $\Delta: X\rightarrow X\times X$, there are two projections $p_1$ and $p_2: X\times X\rightarrow X$ such that $p_i\circ\Delta=id$.  This implies that the source map $\pi_1$ and the target map $\pi_2$ are equal in the derived category in this case, so $X\times^R_{X\times X}X$ becomes a group over $X$ \cite{ACH}.  A similar argument works if the inclusion from $X$ to its formal neighborhood in $S$ splits.

Generally speaking, $N_{X/S}[-1]$ has an $L_\infty$ algebroid structure in the $(\infty,1)$-category of dg quasi-coherent sheaves on $X$.  However, the Lie bracket may not be $O_X$-linear, and it may not satisfy the Jacobi identity when we pass to the derived category of $X$.  Calaque, C\u{a}ld\u{a}raru, and Tu proved that the induced bracket in the derived category is $\mathcal{O}_X$-linear if $i$ splits to first order, and it satisfies the Jacobi identity if $i$ splits to second order \cite{CCT}.  As a consequence $N_{X/S}[-1]$ admits a natural Lie algebra structure in the derived category if $i$ splits to second order.   Later, Calaque and Grivaux showed that $N_{X/S}[-1]$ has a natural Lie algebra structure if $X\hookrightarrow S$ is a tame quantized cycle \cite{CG}, a weaker condition than splitting to second order.  More precisely, for an embedding $i:X\hookrightarrow S$ with a chosen first order splitting, they described the $\mathcal{O}_X$-linear bracket $N_{X/S}\otimes N_{X/S}\rightarrow\S^2N_{X/S}\rightarrow N_{X/S}[1]$ explicitly as the extension class of the short exact sequence of vector bundles on $X$
$$0\rightarrow\S^2N_{X/S}^\vee\cong\frac{I_X^2}{I_X^3}\rightarrow \varphi_*\frac{I_X}{I_X^3}\rightarrow \frac{I_X}{I_X^2}\cong N_{X/S}^\vee\rightarrow0,$$
where $I_X$ is the ideal sheaf of $X$ in $S$.  This bracket satisfies the Jacobi identity under the tameness assumption.  In the rest of the paper we will only use embeddings which are split to first order without requiring the Jacobi identity to hold for this specific bracket.  We will call this type of bracket to be a pre-Lie bracket.

\paragraph We finished background materials.  The rest of this section is about the definition of the Bass-Quillen class and the Lie theoretic interpretations of our results.  From now on we fix smooth subvarieties $X$ and $Y$ of $S$ with closed embeddings $i$, $j$, and $f$
$$
\xymatrix{
X \ar@{^{(}->}[r]^{f}\ar@{^{(}->}@/_1.0pc/[rr]_{i} & Y\ar@{^{(}->}[r]^{j} & S.
}
$$
We assume that $f$, $j$, and $i$ are split to first order, and that we have fixed first order splittings of $f$, $j$, and $i$ which will be denoted as $s$, $\pi$, and $\varphi$.

\paragraph It is crucial to note that all the constructions of HKR isomorphisms depend on the choice of splitting.  Therefore, we need to assume some compatibility on these splittings, namely, that $p\circ q|_X=\rho$
$$
\xymatrix{
T_X \ar[r] & T_Y|_X \ar[r]\ar@{-->}@/_1.0pc/[l]_{p} & T_S|_X\ar@{-->}@/_1.0pc/[l]_{q|_X}\ar@{-->}@/^1.0pc/[ll]^{\rho},
}
$$
where $p$, $q$, and $\rho$ are splittings on the tangent bundles corresponding to the first order splittings $s$, $\pi$, and $\varphi$.

\paragraph{\bf Definition of the Bass-Quillen class.} There is a class $\alpha_{s,N^\vee_{Y/S}|_{X^{(1)}_Y}} \in\mathrm{Ext}^1(N_{X/Y}\otimes N_{Y/S}|_X, N_{Y/S}|_X)$ which plays an important role in what follows.  It has interpretations both in terms of Lie theory and as the obstruction for a general Bass-Quillen theorem to hold.

Consider the short exact sequence of $\mathcal{O}_{X^{(1)}_Y}$-modules
$$0\rightarrow t_*N^\vee_{X/Y}\rightarrow\mathcal{O}_{X^{(1)}_Y}\rightarrow t_*\mathcal{O}_X\rightarrow0,$$
where $t$ is the inclusion $X\hookrightarrow X^{(1)}_Y$.

For any vector bundle $\mathcal{M}$ on $X^{(1)}_Y$, tensor it with this short exact sequence.  Then push-forward the sequence onto $X$ via $s$.  Using the fact that $s\circ t=id$ and the projection formula, we get a short exact sequence of vector bundles on $X$
$$0\rightarrow N^\vee_{X/Y}\otimes\mathcal{M}|_X\rightarrow s_*\mathcal{M}\rightarrow \mathcal{M}|_X\rightarrow0.$$
Dualizing, we get an extension class $\alpha_{s,\mathcal{M}}: N_{X/Y}\otimes\mathcal{M}^\vee|_X\rightarrow\mathcal{M}^\vee|_X[1]$.  We call $\alpha_{s,\mathcal{M}}$ the Bass-Quillen class for the pair $(s,\mathcal{M})$ for the reason explained in~(\ref{subsec:BQ}).  The class $\alpha_{s,\mathcal{M}}$ vanishes if and only if $\mathcal{M}$ is isomorphic to $s^*\mathcal{M}|_X$ \cite{CG}.

Now, choose $\mathcal{M}$ to be $N^\vee_{Y/S}|_{X^{(1)}_Y}$ which gives the class $\alpha_{s,N^\vee_{Y/S}|_{X^{(1)}_Y}}$.  It admits an interpretation as a Lie module structure map.  We explain this interpretation in the rest of this section.
\paragraph {\bf Lie theoretic interpretations of the Bass-Quillen class and Theorem B.} To make our Lie theoretic interpretation clearer, let us assume that all the three derived self-intersections $X\times^R_Y X$, $X\times^R_S X$, and $Y\times^R_SY$ are groups in this section.  However, we will state and explain theorems and propositions later in sections 3 to 6 which only assume existence of first order splittings.  One can check that the natural maps
$$X\times^R_YX\longrightarrow X\times^R_SX\longrightarrow X\times^R_SY=Y\times^R_SY|_X$$
 are maps of groups.  All the shifted normal bundles are Lie algebras under this assumption.  In this section we denote $N_{X/Y}[-1]$, $N_{X/S}[-1]$, and $N_{Y/S}|_X[-1]$ by $\mathfrak{h}$, $\mathfrak{g}$, and $\mathfrak{n}$ respectively.  The functoriality of HKR isomorphisms can be viewed as the functoriality of the exponential maps from Lie algebras to Lie groups.
\paragraph The map $\mathfrak{h}=N_{X/Y}[-1]\hookrightarrow \mathfrak{g}=N_{X/S}[-1]$ preserves the Lie brackets, so we are able to prove the commutativity of the left square in Theorem A with no difficulty.  Moreover, the compatibility of the Lie brackets implies that $\mathfrak{g}$ is an $\mathfrak{h}$-module, and $\mathfrak{h}\hookrightarrow\mathfrak{g}$ is a map of $\mathfrak{h}$-modules.  Therefore, $\mathfrak{n}=\mathfrak{g}/\mathfrak{h}=N_{Y/S}|_X[-1]$ has a natural $\mathfrak{h}$-module structure.  We get a short exact sequence of $\mathfrak{h}$-modules
$$\xymatrix{
 0\ar[r] & \mathfrak{h}=N_{X/Y}[-1] \ar[r] & \mathfrak{g}=N_{X/S}[-1] \ar[r] & \mathfrak{n}=\mathfrak{g}/\mathfrak{h}=N_{Y/S}|_X[-1]\ar[r]& 0.
 }
$$
The $\mathfrak{h}$-module structure on $\mathfrak{g}/\mathfrak{h}$: $N_{X/Y}\otimes N_{Y/S}|_X\rightarrow N_{Y/S}|_X[1]$ is exactly the Bass-Quillen class $\alpha_{s,N^\vee_{Y/S}|_{X^{(1)}_Y}}$.  We will prove this statement in section 5.\label{subsec: Lie for BQ}




\paragraph On the other hand, the map $\mathfrak{g}=N_{X/S}[-1]\rightarrow \mathfrak{n}=N_{Y/S}|_X[-1]$ may not in general preserve the Lie brackets even if we assume that all the derived self-intersections are groups.  This explains the difficulty for proving the functoriality of the exponential maps in the right square of Theorem A.  In section 6 we will show that $\mathfrak{h}=N_{X/Y}[-1]$ acts on its module $\mathfrak{g}/\mathfrak{h}=N_{Y/S}|_X[-1]$ trivially if and only if the Lie brackets are preserved, i.e., $\mathfrak{g}\rightarrow\mathfrak{g}/\mathfrak{h}=\mathfrak{n}$ is a map of Lie algebras.  This is Theorem B.  As a consequence the right square in Theorem A commutes when the Lie brackets are preserved.

\paragraph Let us restate~(\ref{subsec: Lie for BQ}) and Theorem B for Lie algebras in the category of vector spaces.  By abuse of notations, we use the same notations $\mathfrak{h}$ and $\mathfrak{g}$ for Lie algebras in the category of vector spaces.  Consider an injective morphism of Lie algebras $\mathfrak{h}\hookrightarrow\mathfrak{g}$.  The quotient $\mathfrak{n}=\mathfrak{g}/\mathfrak{h}$ is naturally an $\mathfrak{h}$-module, so we get a short exact sequence of $\mathfrak{h}$-modules
$$
\xymatrix{
 0\ar[r] & \mathfrak{h} \ar[r]^{\alpha} & \mathfrak{g} \ar[r]^{\beta~~~~} & \mathfrak{n}=\mathfrak{g}/\mathfrak{h}\ar[r]& 0.
}
$$
There is a way to construct a pre-Lie bracket on $\mathfrak{n}$ once a splitting $\mathfrak{n}\dashrightarrow\mathfrak{g}$ is chosen.  It becomes a Lie bracket under the tameness assumption \cite{CG}, but we do not need this pre-Lie bracket to be a Lie bracket in our paper.  The morphism $\beta$ preserves the pre-Lie brackets if and only if $\mathfrak{h}$ acts trivially on $\mathfrak{n}$.  Triviality of the $\mathfrak{h}$-module structure on $\mathfrak{n}$ implies that $\mathfrak{n}$ is actually a Lie algebra, the map $\beta$ is a Lie algebra morphism, and the diagram
$$
\xymatrix{
G\ar[r]^{\Phi} & N \\
\mathfrak{g}\ar[r]^{\beta=d\Phi}\ar[u]^{\mathrm{exp}} & \mathfrak{n}\ar[u]_{\mathrm{exp}}
}
$$
is commutative.

\section{The proof of Theorem A}

In this section we prove Theorem A.
\paragraph Suppose $i: X\hookrightarrow S$ is a closed embedding of smooth schemes with a first order splitting
$$
\xymatrix{
   X\ar@{^{(}->}[r]_{\mu} & X^{(1)}_S \ar@{^{(}->}[r]_{\nu}\ar@{-->}@/_1.0pc/[l]_{\varphi} & S.
   }
$$
In section 2 we recalled the construction of the HKR isomorphism $i^*i_*\mathcal{O}_X\cong\S(N^\vee_{X/S}[1])$ from \cite{AC}.  It is defined as the composite map
$$\xymatrix{
\mu^*\nu^*\nu_*\mu_*\mathcal{O}_X\ar[r] & \mu^*\mu_*\mathcal{O}_X\ar[r]^{\cong} & \mathrm{T}^c(N^\vee_{X/S}[1])\ar[r]^{\exp} & \mathrm{T}(N^\vee_{X/S}[1])\ar[r] &\S(N^\vee_{X/S}[1]).
}
$$
It is easy to see that all the constructions are canonical except for the isomorphism $\mu^*\mu_*\mathcal{O}_X\cong\mathrm{T}^c(N^\vee_{X/S}[1])$ which depends on the choice of the splitting $\varphi$.

\paragraph We have a commutative diagram
$$
\xymatrix{
X\ar[ddr]^{t}\ar[ddrr]^{f}\ar[dddr]_{\mu} & & & \\
 & & & \\
 & X^{(1)}_Y\ar[r]^{a}\ar[d]^{g} \ar@{-->}@/_3.0pc/[uul]_{s}& Y\ar[d]^{b}& \\
 & X^{(1)}_S\ar[r]^{f^{(1)}} \ar@{-->}@/^2.0pc/[uuul]^{\varphi} \ar@{-->}@/_1.5pc/[u]^{\sigma} & Y^{(1)}_S\ar[r] \ar@{-->}@/_1.5pc/[u]^{\pi} & S
}
$$
under the assumptions in Theorem A.  The solid arrows are the obvious ones.  The dotted arrows $\pi$, $s$, and $\varphi$ are the first order splittings of the closed embeddings $j: Y\rightarrow S$, $f: X\rightarrow Y$, and $i: X\rightarrow S$ respectively.  Notice that $X^{(1)}_Y$ is the fiber product of $X^{(1)}_S$ and $Y$ over $Y^{(1)}_S$, so we can pull $\pi$ back along the morphism $f^{(1)}$ to define $\sigma$.  The compatibility condition on the splittings means that $s\circ\sigma=\varphi$.

We remind the readers of the claim of Theorem A which will be proved below.

The left square below is commutative.  If the Bass-Quillen class $\alpha_{s,N^\vee_{Y/S}|_{X^{(1)}_Y}}$ is zero, then the right square below is also commutative.
$$
\xymatrix{
X\times^R_YX\ar[r] & X\times^R_SX\ar[r] & Y\times^R_SX=Y\times^R_SY|_X  \\
\mathbb{N}_{X/Y}[-1]\ar[r]\ar[u]^{\cong} & \mathbb{N}_{X/S}[-1]\ar[r]\ar[u]^{\cong} & \mathbb{N}_{Y/S}|_X[-1]\ar[u]^{\cong},
}
$$
where the horizontal maps between the normal bundles are linear, i.e., they are vector bundle maps.

\begin{Proof}[Proof of Theorem A.] \label{Proof A (2)}To check the commutativity of the left square in Theorem A, it suffices to show that the diagram
$$
\xymatrix{
 \mu^*\mu_*\mathcal{O}_X\ar[r]^{\cong}\ar[d] & \mathrm{T}^c(N^\vee_{X/S}[1])\ar[d] \\
t^*t_*\mathcal{O}_X\ar[r]^{\cong} &\mathrm{T}^c(N^\vee_{X/Y}[1])
}
$$
is commutative since all the other constructions are canonical.  The right vertical map is obtained from the natural vector bundle map $N^\vee_{X/S}\rightarrow N^\vee_{X/Y}$.  The horizontal isomorphisms are constructed using the splittings, from explicit resolutions of $\mu_*\mathcal{O}_X$ and $t_*\mathcal{O}_X$ on $X^{(1)}_S$ and $X^{(1)}_Y$ respectively.  These resolutions are of the form $(\mathrm{T}^c(\varphi^*N^\vee_{X/S}[1]),d)$ and $(\mathrm{T}^c(s^*N^\vee_{X/Y}[1]),d)$ as explained in~(\ref{subsec: definition of HKR_1}).

We have $g^*\varphi^*N^\vee_{X/S}=s^*N^\vee_{X/S}$ using the fact that $\varphi=s\circ\sigma$ and $\sigma\circ g=id$.  There is a natural map of vector bundles $ g^*\varphi^*N^\vee_{X/S}=s^*N^\vee_{X/S}\rightarrow s^*N^\vee_{X/Y}$ which induces a map of complexes $g^*(\mathrm{T}^c(\varphi^*N^\vee_{X/S}[1]),d)\rightarrow(\mathrm{T}^c(s^*N^\vee_{X/Y}[1]),d)$.  One can check carefully the induced map is indeed a map of complexes, i.e., the differentials are preserved.  This proves that the diagram
$$
\xymatrix{
  g^*(\mathrm{T}^c(\varphi^*N^\vee_{X/S}[1]),d)\ar[r]\ar[d] & g^*\mu_*\mathcal{O}_X \ar[r]\ar[d] & 0\\
(\mathrm{T}^c(s^*N^\vee_{X/Y}[1]),d)\ar[r] & t_*\mathcal{O}_X \ar[r] & 0
}
$$
which relates the two explicit resolutions of $\mathcal{O}_X$ as an $\mathcal{O}_{X^{(1)}_S}$-algebra and as an $\mathcal{O}_{X^{(1)}_Y}$-algebra is commutative.  If we pull the natural map $g^*\varphi^*N^\vee_{X/S}=s^*N^\vee_{X/S}\rightarrow s^*N^\vee_{X/Y}$ back to $X$, we get the natural vector bundle map $N^\vee_{X/S}\rightarrow N^\vee_{X/Y}$.  This proves that we get our desired commutative diagram at the beginning of the proof of Theorem A once we pull the commutative diagram above back to $X$.

\paragraph\label{subsec 3.3} Similarly, to prove the commutativity of the right square of Theorem A, it suffices to show that the diagram
$$
\xymatrix{
 \mu^*\mu_*\mathcal{O}_X\ar[r]^{\cong} & \mathrm{T}^c(N^\vee_{X/S}[1]) \\
f^*b^*b_*\mathcal{O}_Y\ar[r]^{\cong}\ar[u] &f^*\mathrm{T}^c(N^\vee_{Y/S}[1])\ar[u]
}
$$
is commutative.  The right vertical map is induced by the natural map of vector bundles $N^\vee_{Y/S}|_X\rightarrow N^\vee_{X/S}$.

 If the Bass-Quillen class $\alpha_{s,N^\vee_{Y/S}|_{X^{(1)}_Y}}$ vanishes, then we have an isomorphism between $a^*N^\vee_{Y/S}$ and $s^*(N^\vee_{Y/S}|_X)$.  The latter maps to $s^*N^\vee_{X/S}$ naturally.  Therefore, we get a map $\sigma^*a^*N^\vee_{Y/S}\cong\sigma^*s^*N^\vee_{Y/S}|_X\rightarrow\sigma^*s^*N^\vee_{X/S}=\varphi^*N^\vee_{X/S}$.  Notice that $a\circ\sigma=\pi\circ f^{(1)}$ by the definition of $\sigma$, so we get a map $f^{(1)*}\pi^*N^\vee_{Y/S} =\sigma^*a^*N^\vee_{Y/S}\rightarrow\sigma^*s^*N^\vee_{X/S}=\varphi^*N^\vee_{X/S}$.  This map induces a map of complexes $f^{(1)*}(\mathrm{T}^c(\pi^*N^\vee_{Y/S}[1]),d)\rightarrow(\mathrm{T}^c(\varphi^*N^\vee_{X/S}[1]),d)$.  As a consequence the diagram of resolutions of $\mu_*\mathcal{O}_X$ and $b_*\mathcal{O}_Y$
 $$
 \xymatrix{
   (\mathrm{T}^c(\varphi^*N^\vee_{X/S}[1]),d)\ar[r] & \mu_*\mathcal{O}_X \ar[r] & 0\\
f^{(1)*}(\mathrm{T}^c(\pi^*N^\vee_{Y/S}[1]),d)\ar[r]\ar[u] & f^{(1)*}b_*\mathcal{O}_Y \ar[r]\ar[u] & 0
 }
 $$
is commutative.  We recover the map of vector bundles $N^\vee_{Y/S}|_X\rightarrow N^\vee_{X/S}$ if we pull the natural map $\sigma^*a^*N^\vee_{Y/S}=f^{(1)*}\pi^*N^\vee_{Y/S}\rightarrow\varphi^*N^\vee_{X/S}=\sigma^*s^*N^\vee_{X/S}$  back to $X$.  This proves that we get our desired commutative diagram at the beginning of~(\ref{subsec 3.3}) once we pull the commutative diagram above back to $X$.
\end{Proof}

\section{Theorem B in classical Lie theory}

As a warm-up to proving Theorem B, we present here an analogous result in Lie theory.  We give a proof of this result using techniques that can be adapted to the derived setting of Theorem B.

\paragraph Consider an injective map of Lie algebras in vector spaces $\alpha:\mathfrak{h}\hookrightarrow\mathfrak{g}$.  There is a short exact sequence of $\mathfrak{h}$-modules
$$\xymatrix{
0\ar[r] & \mathfrak{h}\ar[r]^{\alpha} & \mathfrak{g}\ar[r]^{\beta~~~~~~} & \mathfrak{n}=\mathfrak{g}/\mathfrak{h}\ar[r] & 0.
}
$$
Given a vector space map $\gamma:\mathfrak{n}\dashrightarrow\mathfrak{g}$ splitting $\beta$ we define a pre-Lie bracket on $\mathfrak{n}$ by the formula $[x,y]_{\mathfrak{n}}=\beta([\gamma(x),\gamma(y)]_{\mathfrak{g}})$ for any $x,y\in\mathfrak{n}$.  In general, the map $\beta$ may not respect the pre-Lie brackets.

We define a map $\mathfrak{g}\otimes\mathfrak{n}\rightarrow\mathfrak{n}$
$$\sum_i x_i\otimes y_i \rightarrow \sum_i\beta([x_i,\gamma(y_i)]),$$
for $x_i\in\mathfrak{g}$ and $y_i\in\mathfrak{n}$.  This map may not define a $\mathfrak{g}$-module structure on $\mathfrak{n}$ if $\beta$ is not a morphism of Lie algebras.  We state a proposition which is important in sections 4 and 6.  Its proof is left to the reader.
\begin{Proposition}\label{prop Lie observation} In general, we do not expect
$$\xymatrix{
\mathfrak{g}\otimes\mathfrak{g}\ar[r]\ar[d]^{id\otimes\beta} & \mathfrak{g}\ar[d]^{\beta}\\
\mathfrak{g}\otimes\mathfrak{n}\ar[r] & \mathfrak{n}
}$$
to be commutative.  However, the diagram
$$\xymatrix{
\wedge^2\mathfrak{g}=\wedge^2\mathfrak{h}\oplus\wedge^2\mathfrak{n}\oplus(\mathfrak{h}\otimes\mathfrak{n}) \ar[r]\ar[d] & \mathfrak{g}\ar[d]^{\beta}\\
(\mathfrak{h}\otimes\mathfrak{n})\oplus\wedge^2\mathfrak{n}\ar[r] & \mathfrak{n}
}$$
is commutative if we identify $\mathfrak{g}$ with $\mathfrak{h}\oplus\mathfrak{n}$ as a direct sum of vector spaces via $\gamma$.

Therefore the right hand side square of the diagram
$$\xymatrix{
\mathfrak{g}\otimes\mathfrak{g}\ar[r]\ar[d]^{id\otimes\beta} & \wedge^2\mathfrak{g}=\wedge^2\mathfrak{h}\oplus\wedge^2\mathfrak{n}\oplus(\mathfrak{h}\otimes\mathfrak{n}) \ar[r]\ar[d] & \mathfrak{g}\ar[d]^{\beta}\\
\mathfrak{g}\otimes\mathfrak{n}\ar[r] &(\mathfrak{h}\otimes\mathfrak{n})\oplus\wedge^2\mathfrak{n}\ar[r] & \mathfrak{n}
}$$
is commutative, but the square on the left is not.
\end{Proposition}\\

Here is the analogous theorem to Theorem B.
\begin{Theorem}\label{Thm classical Lie} The map $\beta$ preserves the pre-Lie brackets of $\mathfrak{g}$ and $\mathfrak{n}$ if and only if $\mathfrak{h}$ acts trivially on $\mathfrak{n}$.  This is also equivalent to saying that $\beta$ is a morphism of Lie algebras.
\end{Theorem}

\begin{Proof} It is easy to see that the map $\mathfrak{g}\otimes\mathfrak{n}\rightarrow\mathfrak{n}$ has the following properties.

$(I)$ It is compatible with the $\mathfrak{h}$-module structure on $\mathfrak{n}$.  Equivalently, the diagram
$$\xymatrix{
\mathfrak{h}\otimes\mathfrak{n}\ar[r]\ar[d]^{\alpha\otimes id} & \mathfrak{n}\ar[d]^{id}\\
\mathfrak{g}\otimes\mathfrak{n}\ar[r] &\mathfrak{n} }
$$
is commutative.  This follows from the observation that the $\mathfrak{h}$-module structure on $\mathfrak{n}$ can be defined as
$$x\otimes y\rightarrow \beta([\alpha(x),\gamma(y)]),$$
for any $x\in\mathfrak{h}$ and $y\in\mathfrak{n}$.

$(II)$ It defines a $\mathfrak{g}$-module structure on $\mathfrak{n}$ if $\beta$ is a morphism of Lie algebras.  Equivalently, the diagram
$$\xymatrix{
\mathfrak{g}\otimes\mathfrak{n}\ar[r]\ar[d]^{\beta\otimes id} &\mathfrak{n}\ar[d]^{id}\\
\mathfrak{n}\otimes\mathfrak{n}\ar[r] & \mathfrak{n}
}
$$
is commutative if $\beta$ is a morphism of Lie algebras.

$(III)$ The diagram
$$\xymatrix{
\mathfrak{n}\otimes\mathfrak{n}\ar[r]\ar@{-->}[d]^{\gamma\otimes id} &\wedge^2\mathfrak{n}\ar[r]\ar[d] & \mathfrak{n}\\
\mathfrak{g}\otimes\mathfrak{n}\ar[r] &(\mathfrak{h}\otimes\mathfrak{n})\oplus\wedge^2\mathfrak{n}\ar[r] & \mathfrak{n}\ar[u]
}
$$
is commutative.

We first prove that $\mathfrak{h}$ acts on $\mathfrak{n}$ trivially if $\beta$ preserves the pre-Lie brackets.  Compose the two commutative diagrams in $(I)$ and $(II)$ together
$$\xymatrix{
\mathfrak{h}\otimes\mathfrak{n}\ar[r]\ar[d]^{\alpha\otimes id} & \mathfrak{n}\ar[d]^{id}\\
\mathfrak{g}\otimes\mathfrak{n}\ar[r]\ar[d]^{\beta\otimes id} &\mathfrak{n}\ar[d]^{id}\\
\mathfrak{n}\otimes\mathfrak{n}\ar[r] & \mathfrak{n}.
}
$$
The map $\mathfrak{h}\otimes\mathfrak{n}\rightarrow\mathfrak{n}$ is zero because $\beta\circ\alpha=0$.

Let us turn to prove that $\beta$ preserves the pre-Lie brackets if the Lie module structure map $\mathfrak{h}\otimes\mathfrak{n}\rightarrow\mathfrak{n}$ is zero.  Identify $\mathfrak{g}$ with $\mathfrak{h}\oplus\mathfrak{n}$ via $\gamma$.  With property $(I)$ and $(III)$ one can conclude that the map $\mathfrak{g}\otimes\mathfrak{n}=(\mathfrak{h}\otimes\mathfrak{n})\oplus(\mathfrak{n}\otimes\mathfrak{n})\rightarrow\mathfrak{n}$ that we defined at the beginning of this section is the Lie module structure map $\mathfrak{h}\otimes\mathfrak{n}\rightarrow\mathfrak{n}$ plus the pre-Lie bracket on $\mathfrak{n}$: $\mathfrak{n}\otimes\mathfrak{n}\rightarrow\mathfrak{n}$.

If $\mathfrak{h}\otimes\mathfrak{n}\rightarrow\mathfrak{n}$ is zero, then the diagram
$$\xymatrix{
\mathfrak{g}\otimes\mathfrak{n}\ar[r]\ar[d] & (\mathfrak{h}\otimes\mathfrak{n})\oplus\wedge^2\mathfrak{n}\ar[r]\ar[d] & \mathfrak{n}\ar[d]\\
\mathfrak{n}\otimes\mathfrak{n}\ar[r] & \wedge^2\mathfrak{n}\ar[r] & \mathfrak{n}
}
$$
is commutative.  Put the commutative diagram in Proposition \ref{prop Lie observation} and the diagram above together
$$\xymatrix{
\wedge^2\mathfrak{g}\ar[r]\ar[d] & \mathfrak{g}\ar[d]\\
(\mathfrak{h}\otimes\mathfrak{n})\oplus\wedge^2\mathfrak{n}\ar[r]\ar[d] & \mathfrak{n}\ar[d]\\
\wedge^2\mathfrak{n}\ar[r] & \mathfrak{n}.
}
$$
We conclude that $\beta$ is a map of Lie algebras.\end{Proof}

\section{The Bass-Quillen class as a Lie module structure map}

In this section we prove the Lie theoretic interpretation of the Bass-Quillen class that we explained in section 2.  We begin by stating the result under the assumption that closed embeddings split to first order only.  Then we provide explanations in Lie theoretic terms.  We turn to the proof at last.

\paragraph \label{Lem split} We have a map $\Psi: N_{Y/S}|_X\dashrightarrow T_S|_X\rightarrow N_{X/S}$ because $Y$ splits to first order in $S$.  As a consequence the following short exact sequence splits
$$\xymatrix{
 0\ar[r] & N_{X/Y} \ar[r] & N_{X/S} \ar[r] & N_{Y/S}|_X \ar[r]\ar@{-->}@/_1.0pc/[l]_{\Psi} & 0.
}
$$
The short exact sequence above shifted by negative one is analogous to the sequence 
$$\xymatrix{
0\ar[r] & \mathfrak{h}\ar[r]^{\alpha} & \mathfrak{g}\ar[r]^{\beta~~~~~~} & \mathfrak{n}=\mathfrak{g}/\mathfrak{h}\ar[r] & 0.
} 
$$
in section 4.

Most of this section will be devoted to constructing a map $\kappa:N_{X/Y}\otimes N_{X/S}\rightarrow N_{X/S}[1]$.  This map will be given by the extension class of an explicit short exact sequence.  Using this map we will prove the following proposition.
\begin{Proposition}
\label{prop BQ}
In the same setting as Theorem A, the vector bundle map $N_{X/Y}[-1]\rightarrow N_{X/S}[-1]$ preserves the pre-Lie brackets constructed by Calaque and Grivaux \cite{CG}.  There exists a map $\kappa: N_{X/Y}\otimes N_{X/S}\rightarrow N_{X/S}[1]$ defined explicitly by the extension class of a short exact sequence.  The diagram
$$
\xymatrix{
N_{X/Y}\otimes N_{Y/S}|_X\ar[r] & N_{Y/S}|_X[1]\\
N_{X/Y}\otimes N_{X/S}\ar[r]^{\kappa}\ar[u]\ar[d] & N_{X/S}[1]\ar[u]\ar[d]\\
N_{X/S}\otimes N_{X/S}\ar[r] & N_{X/S}[1]}
$$
is commutative.  Here all the vertical maps are the obvious maps of vector bundles.  The top horizontal map is the Bass-Quillen class $\alpha_{s,N^\vee_{Y/S}|_{X^{(1)}_Y}}$, and the bottom horizontal map is the pre-Lie bracket.
\end{Proposition}\\

\paragraph To explain in Lie theoretic terms the meaning of Proposition \ref{prop BQ}, assume for simplicity that the embedding $X\hookrightarrow S$ satisfies the additional conditions that make $\mathfrak{g}=N_{X/S}[-1]$ a Lie algebra. Then denote by $\mathfrak{h}=N_{X/Y}[-1]$.  It is a subalgebra of $\mathfrak{g}$ because the map $\mathfrak{h}\hookrightarrow\mathfrak{g}$ preserves the brackets by Proposition \ref{prop BQ}.  The bundle $N_{Y/S}|_X[-1]$ can be identified with $\mathfrak{g}/\mathfrak{h}$.  Then the diagram in Proposition \ref{prop BQ} becomes
$$
\xymatrix{
\mathfrak{h}\otimes \mathfrak{g}/\mathfrak{h} \ar[r] & \mathfrak{g}/\mathfrak{h}\\
\mathfrak{h}\otimes \mathfrak{g}\ar[r]^{\kappa}\ar[u]\ar[d] & \mathfrak{g}\ar[u]\ar[d]\\
\mathfrak{g}\otimes\mathfrak{g}\ar[r] & \mathfrak{g}. }
$$

The commutativity of
$$
\xymatrix{
N_{X/Y}\otimes N_{X/S}\ar[r]^{\kappa}\ar[d] & N_{X/S}[1]\ar[d]  &\mathfrak{h}\otimes\mathfrak{g}\ar[d]\ar[r]^{\kappa} & \mathfrak{g}\ar[d]  \\
N_{X/S}\otimes N_{X/S}\ar[r] & N_{X/S}[1]              &\mathfrak{g}\otimes\mathfrak{g}\ar[r] & \mathfrak{g} }
$$
says that the morphism $\kappa$ is the structure map of the natural $\mathfrak{h}$-module structure on $\mathfrak{g}=N_{X/S}[-1]$, where $\mathfrak{h}=N_{X/Y}[-1]$ is the Lie algebra.

The commutativity of the diagram
$$
\xymatrix{
N_{X/Y}\otimes N_{Y/S}|_X\ar[r] & N_{Y/S}|_X[1]        & \mathfrak{h}\otimes\mathfrak{g}/\mathfrak{h}\ar[r] &  \mathfrak{g}/\mathfrak{h}\\
N_{X/Y}\otimes N_{X/S}\ar[r]^{\kappa}\ar[u] & N_{X/S}[1]\ar[u]  & \mathfrak{h}\otimes\mathfrak{g}\ar[r]^{\kappa}\ar[u] & \mathfrak{g}\ar[u]
}
$$
says that the Bass-Quillen class $\alpha_{s,N^\vee_{Y/S}|_{X^{(1)}_Y}}$ at the top of the diagram is the structure map of the $\mathfrak{h}$-module structure on $\mathfrak{g}/\mathfrak{h}=N_{Y/S}|_X[-1]$.


\paragraph Before we prove Proposition \ref{prop BQ}, we have to define the morphism $\kappa: N_{X/Y}\otimes N_{X/S}\rightarrow N_{X/S}[1]$ that appears in the middle of the diagram in Proposition \ref{prop BQ}.  We hope to describe the morphism $\kappa$ explicitly as the extension class of a short exact sequence.  There is a technical detail we need to deal with.  As we can see, the pre-Lie bracket is defined as $\S^2N_{X/S}\rightarrow N_{X/S}[1]$ instead of $N_{X/S}\otimes N_{X/S}\rightarrow N_{X/S}[1]$.  It is easy to describe the short exact sequence corresponding to $\S^2N_{X/S}\rightarrow N_{X/S}[1]$.  However, it is hard to describe explicitly what short exact sequence the morphism $N_{X/S}\otimes N_{X/S}\rightarrow N_{X/S}[1]$ corresponds to.  The same phenomenon appears when we try to define $N_{X/Y}\otimes N_{X/S}\rightarrow N_{X/S}[1]$.  There is an anti-symmetric part $\wedge^2N_{X/Y}\hookrightarrow N_{X/Y}\otimes N_{X/S}$.  We can only define our desired Lie module structure map $\kappa$ via the extension class of a short exact sequence after we kill this anti-symmetric part.  Lemma \ref{Lem antisym} below describes how to kill the anti-symmetric part of $N_{X/Y}\otimes N_{X/S}\cong (N_{X/Y}\otimes N_{X/Y})\oplus (N_{X/Y}\otimes N_{Y/S}|_X)$ canonically.

\begin{Lemma}
\label{Lem antisym}
The vector bundle $\frac{I_X^2}{I_X^3+I_Y^2}$ on $X$ is isomorphic to $\S^2N^\vee_{X/Y}\oplus (N^\vee_{X/Y}\otimes N^\vee_{Y/S}|_X)$, where $I_X$ and $I_Y$ are the ideal sheaves of $X$ and $Y$ in $S$.
\end{Lemma}

\begin{proof}
The cokernel of $\S^2N^\vee_{Y/S}|_X\hookrightarrow\S^2N^\vee_{X/S}$ is isomorphic to $\S^2N^\vee_{X/Y}\oplus (N^\vee_{X/Y}\otimes N^\vee_{Y/S}|_X)$ using the splitting in~(\ref{Lem split}).

There is a commutative diagram on $X$
$$\xymatrix{
 0\ar[r] & \displaystyle\frac{I^2_Y}{I^2_Y I_X} \ar[r]\ar[d]^{\cong} & \displaystyle\frac{I^2_X}{I_X^3} \ar[r]\ar[d]^{\cong}& \displaystyle\frac{I_X^2}{I_X^3+I_Y^2}  \ar[r] & 0\\
 0 \ar[r]& \S^2N^\vee_{Y/S}|_X \ar[r]& \S^2N^\vee_{X/S} \ar[r]&\S^2N^\vee_{X/Y}\oplus (N^\vee_{X/Y}\otimes N^\vee_{Y/S}|_X) \ar[r] &0.
}
$$
The two vertical maps above are isomorphisms, so we can complete the diagram above as an isomorphism of short exact sequences.  This implies our desired isomorphism.\end{proof}

\begin{Definition}\label{def module}  Define the morphism $\kappa: N_{X/Y}\otimes N_{X/S}\rightarrow N_{X/S}[1]$ as follows
$$N_{X/Y}\otimes N_{X/S}\rightarrow \S^2N_{X/Y}\oplus (N_{X/Y}\otimes N_{Y/S}|_X)\cong(\frac{I_X^2}{I_X^3+I_Y^2})^\vee\rightarrow N_{X/S}[1],$$
where the map on the left is the obvious map under the identification $N_{X/S}\cong N_{X/Y}\oplus N_{Y/S}|_X$ in~(\ref{Lem split}), and the map on the right is given by the extension class of the short exact sequence
$$0\rightarrow \frac{I_X^2}{I_X^3+I_Y^2}\rightarrow\varphi_*\frac{I_X}{I^3_X+I^2_Y}\rightarrow \frac{I_X}{I^2_X}\rightarrow0.$$
\end{Definition}

\paragraph We will focus on the proof of Proposition \ref{prop BQ}.  The result will follow from Lemma \ref{Lem Preservance} and Propositions \ref{prop proof 3(1)} and \ref{prop Proof 3(2)} below.

\begin{Lemma}
$\S^2N^\vee_{X/Y}$ is isomorphic to $\frac{I_X^2}{I_X^3+I_XI_Y}$.
\end{Lemma}

\begin{proof} The ideal sheaf of $X$ in $Y$ is $\frac{I_X}{I_Y}\subset\mathcal{O}_Y=\mathcal{O}_S/I_Y$.  Note that $\frac{I_X^n}{I_Y^n}\neq(\frac{I_X}{I_Y})^n\subset\mathcal{O}_S/I_Y$.  It is easy to show that $(\frac{I_X}{I_Y})^n\cong\frac{I_X^n+I_Y}{I_Y}\subset\mathcal{O}_Y=\mathcal{O}_S/I_Y$.  Therefore,
$$\S^2N^\vee_{X/Y}\cong\frac{(\frac{I_X}{I_Y})^2}{(\frac{I_X}{I_Y})^3}\cong\frac{I_X^2+I_Y}{I_X^3+I_Y}\cong\frac{I_X^2}{I_X^2\cap(I_X^3+I_Y)}. $$
We have $I_X^2\cap(I_X^3+I_Y)=(I_X^2\cap I_X^3)+(I_X^2\cap I_Y)$ because $I^3_X\subset I^2_X$.  The equality $I^2_X\cap I_Y=I_XI_Y$ is due to the injective map below
$$N^\vee_{Y/S}|_X=\frac{I_Y}{I_YI_X}\hookrightarrow N^\vee_{X/S}=\frac{I_X}{I_X^2}.$$\end{proof}

\begin{Lemma} The map of vector bundles $N_{X/Y}[-1]\rightarrow N_{X/S}[-1]$ preserves the pre-Lie brackets.
\label{Lem Preservance}
\end{Lemma}

\begin{proof} One can check that the two short exact sequences
$$\xymatrix{
0\ar[r] & \displaystyle\frac{I_X^2}{I_X^3}\ar[r]\ar[d] & \varphi_*\displaystyle\frac{I_X}{I_X^3}\ar[r]\ar[d] & \displaystyle\frac{I_X}{I_X^2}\ar[r]\ar[d] & 0\\
0\ar[r] & \displaystyle\frac{I^2_X}{I^3_X+I_XI_Y}\ar[r] & s_* \displaystyle\frac{I_X}{I_X^3+I_Y}\ar[r] & \displaystyle\frac{I_X}{I_X^2+I_Y}\ar[r] & 0
}
$$
are compatible.\end{proof}
\paragraph On the other hand, the map $N_{X/S}[-1]\rightarrow N_{Y/S}|_X[-1]$ may not preserve the pre-Lie brackets because there is no map $(\pi_*\frac{I_Y}{I_Y^3})|_X\rightarrow \varphi_*\frac{I_X}{I_X^3}$ generally.

We prove the commutativity of the two diagrams in Proposition \ref{prop BQ}.  It is divided into two propositions below.

\begin{Proposition}\label{prop proof 3(1)}
The map in Definition \ref{def module} is compatible with the pre-Lie bracket of $N_{X/S}[-1]$.  It is equivalent to saying that the diagram
$$\xymatrix{N_{X/Y}\otimes N_{X/S}\ar[d]\ar[r] &  N_{X/S}[1]\ar[d]\\
 N_{X/S}\otimes N_{X/S}\ar[r]  & N_{X/S}[1]
}$$
is commutative.
\end{Proposition}
\begin{proof} We need to show that the three small diagrams
$$\xymatrix{
N_{X/Y}\otimes N_{X/S}\ar[d]\ar[r] &\S^2N_{X/Y}\oplus (N_{X/Y}\ar[d]\otimes N_{Y/S}|_X)\ar[r]^{~~~~~~~~~~\cong} & \displaystyle(\frac{I_X^2}{I_X^3+I_Y^2})^\vee\ar[d]\ar[r] & N_{X/S}[1]\ar[d]\\
 N_{X/S}\otimes N_{X/S}\ar[r] &\S^2N_{X/S}\ar[r]^{\cong} & \displaystyle(\frac{I_X^2}{I_X^3})^\vee\ar[r]  & N_{X/S}[1]
}
$$
are commutative.  Clearly the one on the left is commutative.  The commutativity of the isomorphism in the middle follows from the compatibility of the two short exact sequences in Lemma \ref{Lem antisym}.  The diagram on the right commutes because the two short exact sequences
$$\xymatrix{
0\ar[r] & \displaystyle\frac{I_X^2}{I_X^3+I^2_Y} \ar[r] & \varphi_* \displaystyle\frac{I_X}{I_X^3+I^2_Y}\ar[r] & \displaystyle\frac{I_X}{I_X^2}\ar[r] & 0\\
0\ar[r] & \displaystyle\frac{I_X^2}{I_X^3}\ar[u]\ar[r] & \varphi_* \displaystyle\frac{I_X}{I_X^3}\ar[u]\ar[r] & \displaystyle\frac{I_X}{I_X^2}\ar[r]\ar[u] & 0
}
$$
are compatible.\end{proof}

\begin{Proposition}\label{prop Proof 3(2)}
There is a commutative diagram
$$\xymatrix{
N_{X/Y}\otimes N_{X/S} \ar[r]\ar[d]   & N_{X/Y}\otimes N_{Y/S}|_X \ar[d]\\
N_{X/S}[1] \ar[r]  & N_{Y/S}[1],
}
$$
where the left vertical map is in Definition \ref{def module}, and the right vertical map is the Bass-Quillen class $\alpha_{s,N^\vee_{Y/S}|_{X^{(1)}_Y}}$.
\end{Proposition}
\begin{proof}  We need to prove that the three small diagrams in the diagram $(1)$
$$\xymatrix{
N_{X/Y}\otimes N_{X/S}\ar[d]\ar[r] & N_{X/Y}\otimes N_{Y/S}|_X\ar[d]\\
\S^2N_{X/Y}\oplus (N_{X/Y}\otimes N_{Y/S}|_X)\ar[d]^{\cong}\ar[r] & N_{X/Y}\otimes N_{Y/S}|_X\ar[d]^{\cong} & (1)\\
\displaystyle(\frac{I_X^2}{I_X^3+I_Y^2})^\vee\ar[d]\ar[r] & \displaystyle(\frac{I_X}{I_X^2+I_Y})^\vee\otimes(\frac{I_Y}{I_YI_X})^\vee\ar[d]\\
N_{X/S}[1]\ar[r] & N_{Y/S}[1]
}
$$
are commutative.  Obviously the one on the top is commutative.

To prove the isomorphism in the middle is compatible, we construct a commutative diagram
$$
\xymatrix{
N^\vee_{X/Y}\otimes N^\vee_{Y/S}|_X \ar@{-->}[r]\ar[d]^{\cong} & N^\vee_{X/S}\otimes N^\vee_{X/S}\ar[r]\ar[d]^{\cong} & \S^2N^\vee_{X/S}\ar[r]\ar[d]^{\cong} & \S^2N^\vee_{X/Y}\oplus(N^\vee_{X/Y}\otimes N^\vee_{Y/S}|_X)\ar[d]^{\cong}\\
\displaystyle{\frac{I_X}{I_X^2+I_Y}\otimes\frac{I_Y}{I_YI_X}}\ar@{-->}[r] & \displaystyle\frac{I_X}{I_X^2}\otimes\frac{I_X}{I_X^2}\ar[r] &\displaystyle\frac{I_X^2}{I_X^3}\ar[r] & \displaystyle\frac{I_X^2}{I_X^3+I_Y^2},
}
$$
where the dotted arrows are defined by the splitting in~(\ref{Lem split}), and the right square commutes as mentioned in the proof of Proposition \ref{prop proof 3(1)}.  Clearly, the left and middle squares are commutative.  We hope to prove that this big commutative diagram is exactly dual to the one in the middle of $(1)$.  It suffices to show that the map
$$\xymatrix{
\displaystyle{\frac{I_X}{I_X^2+I_Y}\otimes\frac{I_Y}{I_YI_X}}\ar@{-->}[r] & \displaystyle{\frac{I_X}{I_X^2}\otimes\frac{I_X}{I_X^2}}\ar[r] &\displaystyle\frac{I_X^2}{I_X^3}\ar[r] & \displaystyle\frac{I_X^2}{I_X^3+I_Y^2}
}
$$
defined using the splitting is equal to the natural map
$$\frac{I_X}{I_X^2+I_Y}\otimes\frac{I_Y}{I_YI_X}\rightarrow \frac{I_X^2}{I_X^3+I_Y^2},$$
$$(a\otimes b)\rightarrow ab,\mbox{ }\mbox{for}\mbox{ }a\in\frac{I_X}{I_X^2+I_Y},\mbox{ }\mbox{and}\mbox{ }b\in\frac{I_Y}{I_YI_X}.$$
One can check this easily.

Let us focus on the commutativity of the bottom square in $(1)$
$$\xymatrix{
 \displaystyle(\frac{I_X^2}{I_X^3+I_Y^2})^\vee\ar[r]\ar[d] & N_{X/S}[1]\ar[d]\\
 \displaystyle{(\frac{I_X}{I_X^2+I_Y})^\vee\otimes(\frac{I_Y}{I_YI_X})^\vee}\ar[r] & N_{Y/S}[1].
}
$$
The bottom horizontal map is defined by a short exact sequence
$$0\rightarrow N^\vee_{X/Y}\otimes N^\vee_{Y/S}|_X\rightarrow s_*a^*N^\vee_{Y/S}\rightarrow N^\vee_{Y/S}\rightarrow0.$$
Moreover, we have a morphism of short exact sequences
$$\xymatrix{
0\ar[r] & N^\vee_{X/Y}\otimes N^\vee_{Y/S}|_X\ar[r]\ar[d] & s_*a^*N^\vee_{Y/S}\ar[r]\ar[d]^{\cong} & N^\vee_{Y/S}|_X\ar[r]\ar[d]^{\cong} & 0\\
0\ar[r] & \displaystyle\frac{I_X I_Y}{I_Y(I^2_X+I_Y)}\ar[r] & s_*\displaystyle\frac{I_Y}{I_Y(I_X^2+I_Y)}\ar[r] & \displaystyle\frac{I_Y}{I_YI_X}\ar[r] & 0.
}
$$
This implies that $N^\vee_{X/Y}\otimes N^\vee_{Y/S}\cong\frac{I_X I_Y}{I_Y(I^2_X+I_Y)}$.  It suffices to show that the short exact sequence $0\rightarrow  \frac{I_X I_Y}{I_Y(I^2_X+I_Y)}\rightarrow s_*\frac{I_Y}{I_Y(I_X^2+I_Y)}\rightarrow \frac{I_Y}{I_YI_X}\rightarrow 0$ is compatible with the short exact sequence $0\rightarrow \frac{I_X^2}{I_X^3+I^2_Y} \rightarrow \varphi_* \frac{I_X}{I_X^3+I^2_Y}\rightarrow \frac{I_X}{I_X^2}\rightarrow 0$.

There is a natural map of sheaves $g_*\frac{I_Y}{I_Y(I_X^2+I_Y)}\rightarrow\frac{I_X}{I_X^3+I^2_Y}.$  We get $s_*\frac{I_Y}{I_Y(I_X^2+I_Y)}\rightarrow\varphi_*\frac{I_X}{I_X^3+I^2_Y}$ by applying $\varphi_*$ on both sides, so the two short exact sequences
$$\xymatrix{
0\ar[r] & \displaystyle\frac{I_X I_Y}{I_Y(I^2_X+I_Y)}\ar[r]\ar[d] & s_*\displaystyle\frac{I_Y}{I_Y(I_X^2+I_Y)}\ar[r]\ar[d] & \displaystyle\frac{I_Y}{I_YI_X}\ar[r]\ar[d] & 0.\\
0\ar[r] & \displaystyle\frac{I_X^2}{I_X^3+I^2_Y} \ar[r] & \varphi_* \displaystyle\frac{I_X}{I_X^3+I^2_Y}\ar[r] & \displaystyle\frac{I_X}{I_X^2}\ar[r] & 0
}
$$
are compatible. \end{proof}

\section{The proof of Theorem B}
We generalize the proof in section 4 to prove Theorem B.  We first define a morphism $N_{X/S}\otimes N_{Y/S}|_X\rightarrow N_{Y/S}|_X[1]$ which is analogous to the map $\mathfrak{g}\otimes\mathfrak{n}\rightarrow\mathfrak{n}$ in section 4.  Then we prove similar statements to properties $(I)$ and $(II)$, and Proposition \ref{prop Lie observation} in section 4.  We prove Theorem B at last.

\paragraph The first thing that we need is the map $N_{X/S}\otimes N_{Y/S}|_X\rightarrow N_{Y/S}|_{X}[1]$ which is the analogue of the map $\mathfrak{g}\otimes\mathfrak{n}\rightarrow\mathfrak{n}$ from section 4. We need to deal with the same technical issue that appears in section 5.  Using the splitting in~(\ref{Lem split}) we see that there is an anti-symmetric part $\wedge^2N_{Y/S}|_X$ in $N_{X/S}\otimes N_{Y/S}|_X=(N_{X/Y}\oplus N_{Y/S}|_X)\otimes N_{Y/S}|_X$.  We need to kill this anti-symmetric part canonically.  Lemma \ref{Lem 6.2} and \ref{Lem antisym 2} describe how to do this.

\begin{Lemma}\label{Lem 6.2} $I^2_Y\cap I^2_XI_Y=I^2_YI_X$.
\end{Lemma}

\begin{proof}
We have $I^2_Y\cap I^3_X=I_Y^2I_X$ because the map $$\S^2N^\vee_{Y/S}|_X=\frac{I_Y^2}{I_Y^2I_X}\rightarrow\S^2N^\vee_{X/S}=\frac{I_X^2}{I_X^3}$$
 is injective.  Then we have
$$I_Y^2I_X\subset I_Y^2\cap I_X^2I_Y\subset I^2_Y\cap I^3_X=I_Y^2I_X.$$
\end{proof}

\begin{Lemma} \label{Lem antisym 2} There is an isomorphism of vector bundles $(\frac{I_XI_Y}{I_X^2I_Y})^\vee\cong\S^2N_{Y/S}|_X\oplus(N_{X/Y}\otimes N_{Y/S}|_X)$ on $X$.
\end{Lemma}

\begin{proof} There is a morphism of short exact sequences
$$\xymatrix{
	0\ar[r] & \displaystyle{\frac{I_Y}{I_YI_X}\otimes\frac{I_Y}{I_YI_X}}\ar[r]\ar[d]  & \displaystyle{\frac{I_X}{I^2_X}\otimes\frac{I_Y}{I_YI_X}} \ar[r]\ar[d]^{v} & \displaystyle{\frac{I_X}{I^2_X+I_Y}\otimes\frac{I_Y}{I_YI_X}}\ar[r]\ar[d]^{u}\ar@{-->}@/_2.0pc/[l]^{\tau}& 0\\
	0\ar[r] & \displaystyle\frac{I^2_Y}{I^2_YI_X} \ar[r] &\displaystyle\frac{I_XI_Y}{I^2_XI_Y}\ar[r] & \displaystyle\frac{I_XI_Y}{I_Y(I_X^2+I_Y)}\ar[r] &0.
}
$$
Everything above is clear except for the injectivity of $\frac{I^2_Y}{I^2_YI_X} \rightarrow \frac{I_XI_Y}{I^2_XI_Y}$.  This is due to Lemma \ref{Lem 6.2}.

\label{subsec: antisym 2'} The short exact sequence on the top is the dual of the sequence of the normal bundles tensored with $N^\vee_{Y/S}|_X$, so it splits naturally.  The map $u$ is an isomorphism as mentioned in the proof of Proposition \ref{prop Proof 3(2)}.  One can construct a splitting $v\circ\tau\circ u^{-1}$ for the short exact sequence on the bottom.  Therefore $\frac{I_XI_Y}{I^2_XI_Y}\cong\frac{I^2_Y}{I^2_YI_X}\oplus\frac{I_XI_Y}{I_Y(I_X^2+I_Y)}\cong\S^2N^\vee_{Y/S}|_X\oplus(N^\vee_{X/Y}\otimes N^\vee_{Y/S}|_X)$.  The diagram of two short exact sequences above says that there is a commutative diagram
$$\xymatrix{
		0\ar[r] & N^\vee_{Y/S}|_X\otimes N^\vee_{Y/S}|_X\ar[r]\ar[d]  & N^\vee_{X/S}\otimes N^\vee_{Y/S}|_X \ar[r]\ar[d]^{v} & N^\vee_{X/Y}\otimes N^\vee_{Y/S}|_X\ar[r]\ar[d]^{id}\ar@{-->}@/_2.0pc/[l]^{\tau}& 0\\
		0\ar[r] & \S^2N^\vee_{Y/S}|_X \ar[r]\ar[d]^{\cong} & \S^2N^\vee_{Y/S}|_X\oplus(N^\vee_{Y/S}|_X\otimes N^\vee_{X/Y})\ar[r]\ar[d]^{\cong} & N^\vee_{X/Y}\otimes N^\vee_{Y/S}|_X\ar[r]\ar[d]^{\cong} &0\\
0\ar[r] & \displaystyle\frac{I_Y^2}{I_Y^2I_X} \ar[r] & \displaystyle\frac{I_XI_Y}{I^2_XI_Y}\ar[r] & \displaystyle\frac{I_XI_Y}{(I_X^2+I_Y)I_Y}\ar[r] &0.
}
$$\end{proof}

\begin{Definition}\label{def module 2} Define the morphism $N_{X/S}\otimes N_{Y/S}|_X\rightarrow N_{Y/S}[1]$ which is analogous to $\mathfrak{g}\otimes\mathfrak{n}\rightarrow\mathfrak{n}$ in section 4 as follows
$$N_{X/S}\otimes N_{Y/S}|_X\rightarrow\S^2(N_{Y/S}|_X)\oplus(N_{X/Y}\otimes N_{Y/S}|_X)\cong(\frac{I_XI_Y}{I_X^2I_Y})^\vee\rightarrow N_{Y/S}|_X[1],$$
where the map on the left hand side is the obvious one under the identification $N_{X/S}\cong N_{X/Y}\oplus N_{Y/S}|_X$, and the map on the right hand side is given by the extension class of the short exact sequence
$$0\rightarrow  \frac{I_XI_Y}{I_X^2I_Y}\rightarrow \varphi_*\frac{I_Y}{I_YI_X^2}\rightarrow \frac{I_Y}{I_YI_X}\rightarrow0.$$
\end{Definition}

The following proposition is analogous to property $(I)$ in section 4.

\begin{Proposition} \label{prop 6.6} There is a commutative diagram
$$\xymatrix{
	N_{X/S}\otimes N_{Y/S}|_X\ar[r]   & N_{Y/S}|_X[1]\\
	 N_{X/Y}\otimes N_{Y/S}|_X\ar[r] \ar[u] & N_{Y/S}|_X[1],\ar[u]
}$$
where the horizontal map at top is in Definition \ref{def module 2}, and the horizontal map at bottom is the Bass-Quillen class $\alpha_{s,N^\vee_{Y/S}|_{X^{(1)}_Y}}$.
\end{Proposition}

\begin{proof} The isomorphism $N^\vee_{X/Y}\otimes N^\vee_{Y/S}|_X\cong\frac{I_XI_Y}{(I_X^2+I_Y)I_Y}$ is mentioned in the proof of Proposition \ref{prop Proof 3(2)}.  It suffices to prove that the three small diagrams
$$
\xymatrix{
N_{X/S}\otimes N_{Y/S}|_X\ar[r] & \S^2N_{Y/S}|_X\oplus(N_{X/Y}\otimes N_{Y/S}|_X)\ar[r]^{~~~~~~~~\cong} & \displaystyle(\frac{I_XI_Y}{I_X^2I_Y})^\vee\ar[r] & N_{Y/S}|_X[1]\\
N_{X/Y}\otimes N_{Y/S}|_X\ar[r]\ar[u] & N_{X/Y}\otimes N_{Y/S}|_X\ar[r]^{\cong}\ar[u] &  \displaystyle(\frac{I_XI_Y}{(I_X^2+I_Y)I_Y})^\vee\ar[r]\ar[u] & N_{Y/S}|_X[1]\ar[u]
}
$$
are commutative.  Obviously, the left one is commutative.  Commutativity of the one in the middle is due the the compatibility of short exact sequence in Lemma \ref{Lem antisym 2}.  The rest of our proof is devoted to the commutativity of the diagram on the right.

There is a natural map $\frac{I_Y}{I_YI_X^2}\rightarrow g_*\frac{I_Y}{(I^2_X+I_Y)I_Y}$.  We get $\varphi_*\frac{I_Y}{I_YI_X^2}\rightarrow s_*\frac{I_Y}{(I^2_X+I_Y)I_Y}$ by applying $\varphi_*$ on both sides.  This gives the two compatible short exact sequences
$$\xymatrix{
	0\ar[r] & \displaystyle\frac{I_XI_Y}{I_X^2I_Y}\ar[r]\ar[d]  & \varphi_* \displaystyle\frac{I_Y}{I_YI_X^2}\ar[r]\ar[d] & \displaystyle\frac{I_Y}{I_YI_X}\ar[r]\ar[d]
	& 0\\
	0\ar[r] & \displaystyle\frac{I_XI_Y}{(I_X^2+I_Y)I_Y} \ar[r] & s_*\displaystyle\frac{I_Y}{(I^2_X+I_Y)I_Y}\ar[r] & \displaystyle\frac{I_Y}{I_YI_X}\ar[r] &0,
}
$$
which proves that the diagram on the right is commutative.\end{proof}

The following proposition is similar to Proposition \ref{prop Lie observation}.

\begin{Proposition} \label{prop Lie obs in $6} There is a commutative diagram
$$\xymatrix{
\S^2N_{X/S}\ar[r]\ar[d] &N_{X/S}[1]\ar[d]\\
\S^2N_{Y/S}|_X\oplus(N_{X/Y}\otimes N_{Y/S}|_X)\cong\displaystyle(\frac{I_XI_Y}{I_X^2I_Y})^\vee\ar[r] & N_{Y/S}|_X[1],
}
$$
where the bottom horizontal map is in Definition \ref{def module 2}.

However, we do not expect the following big diagram
$$\xymatrix{
N_{X/S}\otimes N_{X/S}\ar[r]\ar[d] &\S^2N_{X/S}\ar[r]\ar[d] &N_{X/S}[1]\ar[d]\\
N_{X/S}\otimes N_{Y/S}|_X\ar[r] & \S^2N_{Y/S}|_X\oplus(N_{X/Y}\otimes N_{Y/S}|_X)\cong\displaystyle(\frac{I_XI_Y}{I_X^2I_Y})^\vee\ar[r] & N_{Y/S}|_X[1]
}
$$
is commutative generally.  One can check that the left square is not commutative as mentioned in Proposition \ref{prop Lie observation}.
\end{Proposition}
\begin{proof} We need to prove that the two small diagrams in diagram $(2)$
$$\xymatrix{
\S^2N_{X/S}\ar[r]^{\cong}\ar[d] & \displaystyle(\frac{I_X^2}{I_X^3})^\vee\ar[d]\ar[r] &N_{X/S}[1]\ar[d]\\
\S^2N_{Y/S}|_X\oplus(N_{X/Y}\otimes N_{Y/S}|_X)\ar[r]^{~~~~~~~~~~~~~~\cong} &\displaystyle(\frac{I_XI_Y}{I_X^2I_Y})^\vee\ar[r] & N_{Y/S}|_X[1] & (2)
}
$$
are commutative.  The diagram on the right of $(2)$ commutes because the two short exact sequence
$$\xymatrix{
0\ar[r] & \displaystyle\frac{I_XI_Y}{I_X^2I_Y}\ar[r]\ar[d] & \varphi_*\displaystyle\frac{I_Y}{I_YI_X^2}\ar[r]\ar[d] &\displaystyle\frac{I_Y}{I_YI_X} \ar[r]\ar[d] &0\\
0\ar[r] & \displaystyle\frac{I_X^2}{I_X^3}\ar[r] & \varphi_*\displaystyle\frac{I_X}{I_X^3}\ar[r] &\displaystyle\frac{I_X}{I_X^2} \ar[r] &0.
}
$$
are compatible.  Let us prove that the left diagram in $(2)$ commutes.  We construct a diagram below
$$
\xymatrix{
\S^2N^\vee_{Y/S}\oplus(N^\vee_{X/Y}\otimes N^\vee_{Y/S}|_X)\ar[rd]^{\zeta}\ar@{-->}[r]^{~~~~~~\epsilon} & \S^2N^\vee_{Y/S}|_X\ar[r]^{\epsilon'} & \S^2N^\vee_{X/S} \\
  &N^\vee_{X/Y}\otimes N^\vee_{Y/S}|_X\ar@{-->}[r]^{\zeta'} & N^\vee_{X/S}\otimes N^\vee_{X/S}\ar[u]^{\vartheta}\\
\displaystyle\frac{I_XI_Y}{I_X^2I_Y}\ar@{-->}[r]^{\delta}\ar[rd]^{\xi}\ar[uu] & \displaystyle\frac{I_Y^2}{I_Y^2I_X}\ar[r]^{\delta'}\ar@/_2.0pc/[uu] & \displaystyle\frac{I_X^2}{I_X^3}\ar@/_1.0pc/[uu]\\
   &\displaystyle\frac{I_XI_Y}{(I_X^2+I_Y)I_Y}\cong\displaystyle{\frac{I_X}{I_X^2+I_Y}\otimes \frac{I_Y}{I_YI_X}}\ar@{-->}[r]^{~~~~~~\xi'}\ar@/^1.5pc/[uu] & \displaystyle{\frac{I_X}{I_X^2}\otimes \frac{I_X}{I_X^2}},\ar[u]^{\vartheta'}\ar@/_1.0pc/[uu]
}
$$
where all the vertical maps are natural isomorphisms.  The dotted arrows are constructed by splittings in the proof of Lemma \ref{Lem antisym 2}, and the solid arrows are the obvious ones which also appear in the proof of Lemma \ref{Lem antisym 2}.  The two short exact sequences and their splittings in the proof of Lemma \ref{Lem antisym 2} are compatible, so the diagram above is commutative.  Taking the direct sum and direct product of the maps above, we get a commutative diagram
$$\xymatrix{
\S^2N^\vee_{Y/S}\oplus(N^\vee_{X/Y}\otimes N^\vee_{Y/S}|_X)\ar[r]^{\epsilon\times\zeta} & \S^2N^\vee_{Y/S}\oplus(N^\vee_{X/Y}\otimes N^\vee_{Y/S}|_X) \ar[r]^{~~~~~~~~~~~~~~~~\epsilon'\oplus(\vartheta\circ\zeta')} & \S^2N^\vee_{X/S}\\
\displaystyle\frac{I_XI_Y}{I_X^2I_Y}\ar[r]^{\xi\times\delta}\ar[u]^{\cong} & \displaystyle{\frac{I^2_Y}{I^2_YI_X}\oplus\frac{I_XI_Y}{I_X^2+I_Y}}\ar[r]^{(\vartheta'\circ\xi')\oplus\delta'}\ar[u]^{\cong} & \displaystyle\frac{I^2_X}{I_X^3}\ar[u]^{\cong}.
}
$$
We hope to prove that the diagram above is dual to the left square in $(2)$.  This says that we need to prove the following statement $(3)$:

The map we constructed above $((\vartheta'\circ\xi')\oplus\delta')\circ(\xi\times\delta) $ is equal to the natural map $\frac{I_XI_Y}{I_X^2I_Y}\rightarrow\frac{I_X^2}{I_X^3}$.

Consider the following diagram
$$\xymatrix{
  \displaystyle\frac{I_XI_Y}{I_X^2I_Y}\ar@{-->}[r]^{\delta}\ar[rd]^{\xi} & \displaystyle\frac{I_Y^2}{I_Y^2I_X}\ar[r]^{\delta'} & \displaystyle\frac{I_X^2}{I_X^3}\\
   &\displaystyle\frac{I_XI_Y}{(I_X^2+I_Y)I_Y}\cong\displaystyle{\frac{I_X}{I_X^2+I_Y}\otimes \frac{I_Y}{I_YI_X}}\ar@{-->}[r]_{~~~~~~~~~~~~~~\xi'} & \displaystyle{\frac{I_X}{I_X^2}\otimes \frac{I_X}{I_X^2}}\ar[u]^{\vartheta'}\\
 N^\vee_{X/S}\otimes N^\vee_{Y/S}|_X\ar@{-->}[r]^{\theta}\ar[dr]^{\lambda}\ar[uu] & N^\vee_{Y/S}|_X\otimes N^\vee_{Y/S}|_X\ar@/_2.0pc/[uu]\ar[uur]^{\theta'~~} \\
& N^\vee_{X/Y}\otimes N^\vee_{Y/S}|_X,\ar@/^2.0pc/[uu]\ar[uru]_{\lambda'}
}
$$
where the dotted arrows are the splittings in the proof of Lemma \ref{Lem antisym 2}.  The diagram above is commutative because the two short exact sequences and their splittings in the proof of Lemma \ref{Lem antisym 2} are compatible.

Notice that $N^\vee_{X/S}\otimes N^\vee_{Y/S}|_X\rightarrow \frac{I_XI_Y}{I_X^2I_Y}$ is surjective.  To prove the statement $(3)$, it suffices to show that the map $(\theta'\oplus(\vartheta'\circ\lambda'))\circ(\theta\times\lambda): N^\vee_{X/S}\otimes N^\vee_{Y/S}|_X=\frac{I_X}{I_X^2}\otimes\frac{I_Y}{I_YI_X}\rightarrow \S^2N^\vee_{X/S}=\frac{I_X^2}{I_X^3}$ constructed via the splittings is equal to the natural map
$$\frac{I_X}{I_X^2}\otimes\frac{I_Y}{I_YI_X}\rightarrow\frac{I_X^2}{I_X^3},$$
$$(a\otimes b)\rightarrow ab\mbox{, for }a\in\frac{I_X}{I_X^2},\mbox{ and }b\in\frac{I_Y}{I_YI_X}.$$
One can verify it easily.\end{proof}

The following lemma is analogous to property $(II)$ in section 4.

\begin{Lemma}\label{Lem property II} If $N_{X/S}[-1]\rightarrow N_{Y/S}|_X[-1]$ preserves the pre-Lie brackets, then there is a commutative diagram
$$
\xymatrix{
\S^2N_{Y/S}|_X\ar[r] & N_{Y/S}|_X[1]\\
\S^2N_{Y/S}\oplus(N_{X/Y}\otimes N_{Y/S}|_X)\ar[r]\ar[u] & N_{Y/S}|_X[1]\ar[u]^{id},
}
$$
where the bottom horizontal map is in Definition \ref{def module 2}.
\end{Lemma}

\begin{proof} Put what we want to prove into a larger diagram
$$
\xymatrix{
\S^2N_{Y/S}|_X\ar[rr]^{\chi} & & N_{Y/S}|_X[1]\\
&\S^2N_{Y/S}|_X\oplus(N_{X/Y}\otimes N_{Y/S}|_X)\ar[ul]_{\psi}\ar[ur]^{\phi} & \\
\S^2N_{X/S}\ar[ur]^{\iota}\ar[rr]\ar[uu] & & N_{Y/S}|_X[1].\ar[uu]
}
$$
The outer square commutes because we assume that the brackets are preserved.  We want to show that $\chi\circ\psi=\phi$.  The commutativity of the outer square and Proposition \ref{prop Lie obs in $6} show that $\chi\circ\psi\circ\iota=\phi\circ\iota$.  The map $\iota$ splits naturally, so we have our desired result. \end{proof}

There should be a statement analogous to property $(III)$.  It will appear in the proof of Theorem B below.

\begin{Proof}[Proof of Theorem B.] We first prove that the Bass-Quillen Lie module map $\alpha_{s,N^\vee_{Y/S}|_{X^{(1)}_Y}}$ is zero if the pre-Lie brackets are preserved.  There is a commutative diagram
$$
\xymatrix{
N_{Y/S}|X\otimes N_{Y/S}|_X\ar[r] & \S^2N_{Y/S}|_X\ar[r] & N_{Y/S}|_X[1] \\
 & \S^2N_{Y/S}|_X\oplus(N_{X/Y}\otimes N_{Y/S}|_X)\ar[r]\ar[u] & N_{Y/S}|_X [1]\ar[u]^{id}\\
 & N_{X/S}\otimes N_{Y/S}|_X\ar[r]\ar[u]\ar[uul] & N_{Y/S}|_X[1]\ar[u]^{id} \\
 & N_{X/Y}\otimes N_{Y/S}|_X\ar[r]\ar[u] & N_{Y/S}|_X[1]\ar[u]^{id}
}
$$
due to Lemma \ref{Lem property II} and Proposition \ref{prop 6.6}.  The Bass-Quillen class $\alpha_{s,N^\vee_{Y/S}|_{X^{(1)}_Y}}$ that appears at the bottom of the diagram above vanishes because the composition $N_{X/Y}\rightarrow N_{X/S}\rightarrow N_{Y/S}|_X$ is zero.

Let us turn to prove that the vector bundle map $N_{X/S}[-1]\rightarrow N_{Y/S}|_X[-1]$ preserves the pre-Lie brackets if the Bass-Quillen class $\alpha_{s,N^\vee_{Y/S}|_{X^{(1)}_Y}}$ is zero.  The proof is similar to the one in section 4.  In \cite{CG} Calaque and Grivaux showed that the pre-Lie brackets on $N_{X/S}[-1]$ and $N_{Y/S}|_X[-1]$ defined by the extension classes of short exact sequences can be also defined as follows
$$\xymatrix{
 \S^2N_{X/S}\ar@{-->}[r] & \S^2T_S|_X\ar[r] & T_S|_X[1]\ar[r] & N_{X/S}[1],
}
$$
$$
\xymatrix{
\S^2N_{Y/S}|_X\ar@{-->}[r] & \S^2T_S|_X\ar[r] & T_S|_X[1]\ar[r] & N_{Y/S}|_X[1].
}
$$
The dotted arrow is due the fact that $f:X\hookrightarrow S$ and $j:Y\hookrightarrow S$ split to first order.  The map in the middle is the Atiyah class.

Using the compatibility condition on splittings of tangent bundles and the fact above, we conclude that the diagram
$$\xymatrix{
 \S^2N_{Y/S}|_X\ar@{-->}[r]\ar@{-->}[d] & \S^2T_S|_X\ar[r]\ar[d]^{id} & T_S|_X[1]\ar[r]\ar[d]^{id} & N_{Y/S}|_X[1]\\
\S^2N_{X/S}\ar@{-->}[r] & \S^2T_S|_X\ar[r] & T_S|_X[1]\ar[r] & N_{X/S}[1]\ar[u]
}
$$
is  commutative.  The diagram above and Proposition \ref{prop Lie obs in $6} say that we have a commutative diagram
$$\xymatrix{
 \S^2N_{Y/S}|_X\ar[r]\ar[d] &  N_{Y/S}|_X[1]\ar[d]^{id}\\
 \S^2N_{Y/S}|_X\oplus(N_{Y/S}|_X\otimes N_{X/Y})\ar[r] &  N_{Y/S}|_X[1]
}
$$
which is analogous to property $(III)$ in section 4.  Based on the diagram above and Proposition \ref{prop 6.6} it is clear that the diagram
$$\xymatrix{
 \S^2N_{Y/S}|_X\oplus(N_{Y/S}|_X\otimes N_{X/Y})\ar[r]\ar[d] &  N_{Y/S}|_X[1]\ar[d]^{id}\\
 \S^2N_{Y/S}|_X\ar[r] &  N_{Y/S}|_X[1]
}
$$
is commutative if the Bass-Quillen class $\alpha_{s,N^\vee_{Y/S}|_{X^{(1)}_Y}}$ is zero.  Compose the diagram above with the one in Proposition \ref{prop Lie obs in $6}.  We get a commutative diagram
$$\xymatrix{
 \S^2N_{X/S}\ar[r]\ar[d] &  N_{X/S}[1]\ar[d]\\
 \S^2N_{Y/S}|_X\oplus(N_{Y/S}|_X\otimes N_{X/Y})\ar[r]\ar[d] &  N_{Y/S}|_X[1]\ar[d]^{id}\\
 \S^2N_{Y/S}|_X\ar[r] &  N_{Y/S}|_X[1],
}
$$
so we conclude that the pre-Lie brackets are preserved if the Bass-Quillen class $\alpha_{s,N^\vee_{Y/S}|_{X^{(1)}_Y}}$ is zero.
\end{Proof}

\paragraph We end this section by providing an example where the Bass-Quillen class is not zero.  Consider the embeddings
$$\xymatrix{
X\ar[r]^{\Delta_X~~~~} & X\times X= Y\ar[r]^{\Delta_{X\times X}~~~~}\ar@{-->}@/_1.5pc/[l]_{p_1} & S=X\times X\times X\times X\ar@{-->}@/_1.5pc/[l]_{\pi_1}
}
$$
for a smooth scheme $X$, where all the inclusions are the diagonal embeddings, and all the splittings are the projections to the first factor.  The normal bundle $N_{Y/S}$ is $p_1^*T_X\oplus p_2^* T_X$ in this case, where $p_1$ and $p_2$ are the two projections: $X\times X\rightarrow X$.  It is clear that the Bass-Quillen class $\alpha_{p_1,N^\vee_{Y/S}|_{X^{(1)}_Y}}$: $T_X\otimes(T_X\oplus T_X)\rightarrow T_X[1]$ is zero plus the Atiyah class: $T_X[-1]\otimes T_X[-1]\rightarrow T_X[-1]$.  It is not zero in general.


\begin{thebibliography}{GP}
	
	

\bibitem[AC12]{AC} D. Arinkin, and A. C\u{a}ld\u{a}raru, {\em When is the self-intersection a fibration?}, Adv. Math. $\mathbf{231}$ (2012), no.2, 815-842.
	
\bibitem[ACH19]{ACH} D. Arinkin, A. C\u{a}ld\u{a}raru, and M. Hablicsek, {\em Formality of derived intersections and the orbifold HKR isomorphism}, J. Algebra $\mathbf{540}$ (2019), 100-120.
	
	
\bibitem[CCT14]{CCT} D. Calaque, A. C\u{a}ld\u{a}raru, and J. Tu, {\em On the Lie algebroid of a derived self-intersection}, Adv. Math. $\mathbf{262}$ (2014), 751-783.
	
\bibitem[CG17]{CG} D. Calaque, and J. Grivaux, {\em The Ext algebra of a quantized cycle}, Journal de l'Ecole Polytechnique-Mathematiques $\mathbf{6}$ (2017).
	
	
\bibitem[CR11]{CR} D. Calaque, and C. Rossi, {\em Lectures on Duflo isomorphisms in Lie algebra and complex geometry}, EMS Series of Lectures in Mathematics, European Mathematical Society (2011).
	
\bibitem[Gri14]{Gri} J. Grivaux, {\em Formality of derived intersections}, Doc. Math. $\mathbf{19}$ (2014), 1003-1016.
	
\bibitem[G64]{G} A. Grothendieck, {\em \'El\'ements de g\'eom\'etrie alg\'ebrique. IV. $\acute{E}$tude locale des sch\'emas et des morphismes de sch\'emas.} I., Inst. Hautes \'Etudes Sci. Publ. Math. No. 20 (1964).
	
\bibitem[L81]{L} H. Lindel, {\em On the Bass-Quillen conjecture concerning projective modules over polynomial rings}, Invent. Math. $\mathbf{65}$ (1981), 319-323.
    	
\bibitem[K99]{K} M. Kapranov, {\em Rozansky-Witten weight invariants via Atiyah classes}, Compositio Math. $\mathbf{115}$ (1999), no.1, 71-113.

\bibitem[S96]{S} R.G. Swan, {\em Hochschild cohomology of quasiprojective schemes}, J. Pure Appl. Algebra $\mathbf{110}$ (1996), no.1, 57-80.

\end{thebibliography}
\end{document}